\documentclass[a4paper,12pt]{article}
\usepackage{amssymb}
\usepackage{color}

\newcommand{\R}{\mathbb{R}}

\newcommand{\beq}{\begin{equation} }
\newcommand{\eqq}{\end{equation} }
\newcommand{\cuad}{{\sqcap\kern-.68em\sqcup}}

\newtheorem{teo}{Theorem}[section]

\newtheorem{proposition}{Proposition}[section]

\newtheorem{lemma}{Lemma}[section]
\newtheorem{corollary}{Corollary}[section]
\newtheorem{remark}{Remark}[section]
\newcommand{\bremark}{\begin{remark} \em}
\newcommand{\eremark}{\end{remark} }

\def\beeq{\begin{equation}}
\def\eeq{\end{equation}}
\newcommand{\begeqaet}{\begin{eqnarray*}}
\newcommand{\eneqaet}{\end{eqnarray*}}

\newcommand{\ve}{\varepsilon}

\newcommand{\Om}{{\Omega}}

\begin{document}
\begin{center}{\bf    Existence and uniqueness of positive solutions for a class of logistic type elliptic equations in $\R^N$ involving  fractional Laplacian}\medskip
%%%%%%%%%%%%%%%%%%%%%%%%%%%%%%%%%%%%%%%%%%%%%%%%%%%%%%%%%%%%%%%%%%%%%%
%%%%%%%%%%%%%%%%%%%%%%%%%%%%%%%%%%%%%%%%%%%%%%%%%%%%%%%%%%%%%%%%%%%%%%

\bigskip

\bigskip

{Alexander Quaas and Aliang Xia}

Departamento de  Matem\'atica,  Universidad T\'ecnica Federico Santa Mar\'{i}a

Casilla: V-110, Avda. Espa\~na 1680, Valpara\'{\i}so, Chile.

 {\sl (alexander.quaas@usm.cl and aliangxia@gmail.com)}
\end{center}
%%%%%%%%%%%%%%%%%%%%%%%%%%%%%%%%%%%%%%%%%%%%%%%%%%%%%%%%%%%%%%%%%%%%%%
%%%%%%%%%%%%%%%%%%%%%%%%%%%%%%%%%%%%%%%%%%%%%%%%%%%%%%%%%%%%%%%%%%%%%%
%\keywords{}
%\subjclass{}

\bigskip

%%%%%%%%%%%%%%%%%%%%%%%%%%%%%%%%%%%%%%%%%%%%%%%%%%%%%%%%%%%%%%%%%%%%%%
\begin{abstract}
 In this paper, we study the existence and uniqueness of positive solutions for the following nonlinear fractional elliptic equation:
 \begin{eqnarray*}
(-\Delta)^\alpha u=\lambda a(x)u-b(x)u^p&{\rm in}\,\,\R^N,
 \end{eqnarray*}
 where $ \alpha\in(0,1) $, $ N\ge 2 $, $\lambda >0$, $a$ and $b$ are positive smooth function in $\R^N$ satisfying
 \[
 a(x)\rightarrow a^\infty>0\quad {\rm and}\quad  b(x)\rightarrow b^\infty>0\quad{\rm as}\,\,|x|\rightarrow\infty.
 \]
 Our proof  is based on a comparison principle and existence, uniqueness and asymptotic behaviors  of various 
 boundary blow-up solutions for a class of elliptic equations involving the fractional Laplacian.
\end{abstract}
%%%%%%%%%%%%%%%%%%%%%%%%%%%%%%%%%%%%%%%%%%%%%%%%%%%%%%%%%%%%%%%%%%%%%%
\date{}%\maketitle
%%%%%%%%%%%%%%%%%%%%%%%%%%%%%%%%%%%%%%%%%%%%%%%%%%%%%%%%%%%%%%%%%%%%%%
%%%%%%%%%%%%%%%%%%%%%%%%%%%%%%%%%%%%%%%%%%%%%%%%%%%%%%%%%%%%%%%%%%%%%%

\renewcommand{\thefootnote}{}

\setcounter{equation}{0}
\section{ Introduction}

A celebrated result of Du and Ma \cite{DM} asserts that the uniqueness positive solution of 
\[
-\Delta u=\lambda u-u^p\quad {\rm in}\,\,\R^N
\]
for $N\ge 1$, $\lambda>0$ and $p>1$, is $u\equiv \lambda^{\frac{1}{p-1}}$. Moreover, in \cite{DM}, the authors also consider the following  logistic type equation:
\beq\label{a1}
-\Delta u=\lambda a(x)u-b(x)u^p\quad {\rm in}\,\,\R^N,
\eqq
where $p>1$, $a$ and $b$ are positive smooth function in $\R^N$ satisfying
 \[
 a(x)\rightarrow a^\infty>0\quad {\rm and}\quad  b(x)\rightarrow b^\infty>0\quad{\rm as}\,\,|x|\rightarrow\infty.
 \]
Then they proved that problem (\ref{a1}) has a unique positive solution for each $\lambda>0$. A similar problem for quasi-linear operator has been studied by Du and Guo \cite{DG1}.

In the present  work, we are interested in understanding whether similar results hold for equations involving  a nonlocal diffusion operator, the simplest of which is perhaps 
the fractional Laplacian.  For $\alpha\in (0,1)$, we study the following fractional elliptic problem:
 \begin{eqnarray}\label{1.1}
(-\Delta)^\alpha u=\lambda u-u^p&{\rm in}\,\,\R^N,
 \end{eqnarray}
where $N\ge2$. The fractional Laplacian is defined, up to a normalization constant, by
\[
(-\Delta)^\alpha u(x)=\int_{\R^N}\frac{2u(x)-u(x+y)-u(x-y)}{|y|^{N+2\alpha}}dy,\quad \forall\,\, x\in\R^N.
\]

Our first main result is
\begin{teo}\label{t1}
Let $\lambda >0$.
Suppose 
$u\in C_{loc}^{2\alpha+\beta}(\R^N)\cap L^1(\R^N,\omega)$ for some $\beta>0$ and $\omega=1/(1+|y|^{N+2\alpha})$ is a nonnegative solution of (\ref{1.1}). Then $u$ must be a constant if $p$ verifies
\beq\label{p}
1+2\alpha<p<\frac{1+\alpha}{1-\alpha}.
\eqq
\end{teo}

\begin{remark}
We notice that 
\[
\frac{N+2\alpha}{N-2\alpha}\le\frac{1+\alpha}{1-\alpha},
\]
if $N\ge 2$.
\end{remark}

As in \cite{DM} and \cite{DG1}, our proof of this result based on a comparison principle for concave sublinear problems (see Lemma \ref{l2.1}) and involves boundary blow-up solutions.
We use a rather intuitive squeezing method to proof Theorem \ref{t1} as follows. Denote $B_R$ as a ball centered at the origin with radius $R$. Then problem 
 \begin{eqnarray*}
\left\{\begin{array}{l@{\quad }l}
(-\Delta)^\alpha v=\lambda v-v^p&{\rm in}\,\,B_R,\\
v=0&{\rm in}\,\,\R^N\setminus B_R,
 \end{array}
 \right.
 \end{eqnarray*}
has a unique positive solution $v_R$ if $R$ is large enough for any fixed $\lambda>0$. On the other hand,
the following boundary blow-up propblem
 \begin{eqnarray}\label{b}
\left\{\begin{array}{l@{\quad }l}
(-\Delta)^\alpha w=\lambda w-w^p&{\rm in}\,\,B_R,\\
\lim_{x\in B_R,x\rightarrow\partial B_R}w(x)=+\infty,\\
w=g&{\rm in}\,\,\R^N\setminus \bar{B}_R,
 \end{array}
 \right.
 \end{eqnarray}
for some $g\in L^1(\R^N\setminus \bar{B}_R, \omega)$ and $\lambda>0$, has a  positive solution $w_R$ for any $R>0$. The comparison principle implies
that any entire  positive solution of (\ref{1.1}) satisfies $v_R\le u\le w_R$ in $B_R$.
Moreover, one can show  (see Lemmas \ref{l2.2} and \ref{l2.3} in Section 2) that  both $v_R$ and $w_R$ converge locally uniformly to 
$\lambda^{\frac{1}{p-1}}$ as $R\rightarrow+\infty$. Therefore, $u\equiv \lambda^{\frac{1}{p-1}}$ in $\R^N$. 

Next, we make use of Theorem \ref{t1} to study logistic type fractional elliptic problems with
variable coefficients that are asymptotically positive constants. More precisely, we study the following problem
 \begin{eqnarray}\label{1.2}
(-\Delta)^\alpha u=\lambda a(x)u-b(x)u^p&{\rm in}\,\,\R^N,
 \end{eqnarray}
 where $a$ and $b$ are positive smooth function in $\R^N$. Moreover, we suppose that 
 \beq\label{1.3}
 a(x)\rightarrow a^\infty>0\quad {\rm and}\quad  b(x)\rightarrow b^\infty>0\quad{\rm as}\,\,|x|\rightarrow\infty.
 \eqq
 We can prove that
 
\begin{teo}\label{t2}
Let $\lambda>0$. Suppose $a$ and $b$ are positive smooth function in $\R^N$ and satisfying (\ref{1.3}). Then equation (\ref{1.2})  has a unique positive solution if $p$ verifies (\ref{p}).
\end{teo}
 
We prove Theorem \ref{t2} by a  similar argument as in the proof of Theorem \ref{t2}, we consider the Dirichlet problem and the boundary blow-up problem in a ball $B_R$.
When $R$ is large, these problems have positive solutions $v_R$ and $w_R$ respectively. By comparison principle, as $R\rightarrow \infty$, $v_R$ increase to a minimal 
positive solution of (\ref{1.2}) and $w_R$ decrease to a maximal positive solution of (\ref{1.2}). Therefore, when (\ref{1.2}) has a unique positive solution, $v_R$ and $w_R$ approximate
this unique solution from below and above, respectively.

We mentioned that, in \cite{DM} and \cite{DG1},  the existence and uniqueness results hold provided $p>1$, but in our Theorems \ref{t1} and \ref{t2} we require $p$ satisfying (\ref{p}).
This is because we will use Perron's method (we refer the reader to  User's guide \cite{CIL} for the presentation of Perron's method which extends to the case of nonlocal equations, see for example \cite{BCI, CFQ,FQ})  to construct solution of problem \ref{b} by applying Proposition \ref{p2.2} and choosing 
\[
\tau=-\frac{2\alpha}{p-1}\in (-1,\tau_0(\alpha))
\]
in $V_\tau(x)$ (see (\ref{2.5})). This implies
\[
p<1-\frac{2\alpha}{\tau_0(\alpha)}.
\]
Moreover, in \cite{CHW}, the authors proved that $\tau_0(\alpha)$ has a simplicity formula, that is, $\tau_0(\alpha)=\alpha-1$. Thus, we have
\[
p<1-\frac{2\alpha}{\tau_0(\alpha)}=\frac{1+\alpha}{1-\alpha}.
\]

This article is organized as follows.  In Section 2 we present some preliminary lemmas to prove  a comparison principle
involving the fractional Laplacian, existence and asymptotic behaviors of boundary blow-up solutions.
Section 3 is devoted to prove the existence and uniqueness results of problems (\ref{1.1}) and (\ref{1.2}), i.e., Theorems \ref{t1} and \ref{t2}.

\setcounter{equation}{0}
\section{ Preliminary lemmas}

In this section, we introduce some lemmas which are useful in the proof of our main results. 
The first important ingredient is the comparison principle involving the fractional Laplacian which is useful in dealing with boundary blow-up
problems.

\begin{lemma}\label{l2.1} (Comparison principle) Suppose that $\Omega$ is a bounded domain in $\R^N$, $a(x)$ and $b(x)$ are continuous functions
in $\Omega$ with $\|a\|_{L^\infty(\Omega)}<\infty$ and $b(x)$ nonnegative and not identity zero. Suppose $u_1,u_2\in C^{2\alpha+\beta}(\Omega)$ for some $\beta>0$
are positive in $\Omega$ and satisfy
\beq\label{2.1}
(-\Delta)^\alpha u_1-a(x)u_1+b(x)u_1^p\ge0\ge (-\Delta)^\alpha u_2-a(x)u_2+b(x)u_2^p\quad in\,\, \Omega
\eqq
and $\limsup_{x\rightarrow\partial\Om} (u_2-u_1)\le 0$ with $u_2-u_1\le0$ in $\R^N\setminus \bar{\Omega} $, where $p>1$.
Then $u_2\le u_1$ in $\Omega$.
\end{lemma}

In order to prove Lemma \ref{l2.1}, we need the following proposition.

\begin{proposition}\label{p1}
For $u\ge0$ and $v>0$, we  have
\[
L(u,v)\ge0\quad in\,\,\R^N\times \R^N,
\]
where
\[L(u,v)(x,y)=\left(u(x)-u(y)\right)^2-\left(v(y)-v(x)\right)\left(\frac{u(y)^2}{v(y)}-\frac{u(x)^2}{v(x)}\right).	\]
Moreover, the equality holds if and only if $u=kv$ a.e. for some contant $k$.
\end{proposition}

We note that Proposition \ref{p1} is a special case ($p=2$) of Lemma 4.6 in \cite{PQ} and we omit the proof here. \\

{\bf Proof of Lemma \ref{l2.1}.} 
 Let  $\phi_1 $ and $\phi_2$ be nonnegative functions in $C_0^\infty(\Omega)$. 
By (\ref{2.1}), we obtain that
\begin{eqnarray}\nonumber
\int_{\R^{2N}}\frac{\left(u_1(x)-u_1(y)\right)\left(\phi_1(x)-\phi_1(y)\right)}{|x-y|^{N+2\alpha}}-\frac{\left(u_2(x)-u_2(y)\right)\left(\phi_2(x)-\phi_2(y)\right)}{|x-y|^{N+2\alpha}}dxdy\\\label{2.2}
\ge \int_\Omega b(x)[u_2^p\phi_2-u_1^p\phi_1]dx+\int_\Omega a(x)(u_1\phi_1-u_2\phi_2)dx.
\end{eqnarray}
For $\ve>0$, we denote $\ve_1=\ve$ and $\ve_2=\ve/2$ and let
\[
v_i=\frac{[(u_2+\ve_2)^2-(u_1+\ve_1)^2]^+}{u_i+\ve_i},\quad i=1,2.
\]
By our our assumption, $v_i$ is zero near $\partial \Om$ and in $\R^N\setminus \bar{\Omega}$. Hence $v_i\in X_0^\alpha(D_0)$, where $D_0\subset\subset\Omega$ and $X_0^\alpha(D_0)=\{w\in H^\alpha(\R^N): w=0 \,\, a.e \,\, {\rm in} \,\, \R^N\setminus D_0\}$. In fact,  it is clear that $\|v_1\|_{L^2(\R^N)}=\|v_1\|_{L^2(D_0)}\le C$ and thus it remains to verify that the Gagliardo norm of $v_1$ in $\R^N$ is bounded by a constant. Using the symmetry of the integral in the Gagliardo norm with respect to $x$ and $y$ and the fact that $v_1=0$ in $\R^N\setminus\Om$, we can split as follows
\begin{eqnarray}\nonumber
\int_{\R^N}\int_{\R^N}\frac{|v_1(x)-v_1(y)|^2}{|x-y|^{N+2\alpha}}dxdy&=&\int_{\Om}\int_{\Om}\frac{|v_1(x)-v_1(y)|^2}{|x-y|^{N+2\alpha}}dxdy\\\label{ee}
&+&2\int_{\Om}\left(\int_{\R^N\setminus\Om}\frac{|v_1(x)|^2}{|x-y|^{N+2\alpha}}dy\right)dx.
\end{eqnarray}
Next, we estimate both integrals in the right hand side of  (\ref{ee}) is finite. We first notice that, for any $y\in \R^N\setminus D_0$,
\[
\frac{|v_1(x)|^2}{|x-y|^{N+2\alpha}}=\frac{\chi_{D_0}(x)|v_1(x)|^2}{|x-y|^{N+2\alpha}}\le \chi_{D_0}(x)|v_1(x)|^2\sup_{x\in D_0}\frac{1}{|x-y|^{N+2\alpha}}.
\]
This implies that
\[
\int_{\Om}\left(\int_{\R^N\setminus\Om}\frac{|v_1(x)|^2}{|x-y|^{N+2\alpha}}dy\right)dx\le \left(\int_{\R^N\setminus\Om}\frac{1}{dist(y,\partial D_0)^{N+2\alpha}}dy\right)\|v_1\|_{L^{2}(D_0)}^2<+\infty
\]
since $dist(\partial\Om,\partial D_0)\ge \gamma>0$ and $N+2\alpha>N$. Hence, the second term  in the right hand side of  (\ref{ee}) is finite by the above inequality. 
In order to show the  first term  in the right hand side of  (\ref{ee}) is also finite, we need the following estimates
\begin{eqnarray}\nonumber
\left|\frac{(u_2(x)+\ve_2)^2-(u_1(x)+\ve_1)^2}{u_1(x)+\ve_1}-\frac{(u_2(y)+\ve_2)^2-(u_1(y)+\ve_1)^2}{u_1(y)+\ve_1}\right|\\\nonumber
=\left|\frac{(u_2(x)+\ve_2)^2}{u_1(x)+\ve_1}-\frac{(u_2(y)+\ve_2)^2}{u_1(y)+\ve_1}+(u_1(y)-u_1(x))\right|\\\label{aa}
\le\left|\frac{(u_2(x)+\ve_2)^2}{u_1(x)+\ve_1}-\frac{(u_2(y)+\ve_2)^2}{u_1(y)+\ve_1}\right|+|u_1(y)-u_1(x)|
\end{eqnarray}
and 
\begin{eqnarray}\nonumber
\left|\frac{(u_2(x)+\ve_2)^2}{u_1(x)+\ve_1}-\frac{(u_2(y)+\ve_2)^2}{u_1(y)+\ve_1}\right|\\\nonumber
=\left|\frac{(u_2(x)+\ve_2)^2-(u_2(y)+\ve_2)^2}{u_1(x)+\ve_1}+\frac{(u_2(y)+\ve_2)^2(u_1(y)-u_1(x))}{(u_1(x)+\ve_1)(u_1(y)+\ve_1)}\right|\\\label{aaa}
\le\frac{u_2(x)+u_2(y)+2\ve_2}{u_1(x)+\ve_1}|u_2(x)-u_2(y)|+\frac{(u_2(x)+\ve_2)^2}{(u_1(x)+\ve_1)(u_1(y)+\ve_1)}|u_1(y)-u_1(x)|.
\end{eqnarray}
Combining (\ref{aa}) and (\ref{aaa}), we have
\begin{eqnarray*}
\left|\frac{(u_2(x)+\ve_2)^2-(u_1(x)+\ve_1)^2}{u_1(x)+\ve_1}-\frac{(u_2(y)+\ve_2)^2-(u_1(y)+\ve_1)^2}{u_1(y)+\ve_1}\right|\\
\le C(\ve_1,\ve_2,\|u_1\|_{L^\infty(\Om)},\|u_2\|_{L^\infty(\Om)})(|u_1(y)-u_1(x)|-|u_2(y)-u_2(x)|)\\
\le \tilde{C}|x-y|^{2\alpha+\beta}.
\end{eqnarray*}
In the last inequality of above estimate, we have used the fact $u_1,u_2\in C^{2\alpha+\beta}(\Om)$.
This implies 
\[
\int_{\Om}\int_{\Om}\frac{|v_1(x)-v_1(y)|^2}{|x-y|^{N+2\alpha}}dxdy<+\infty
\]
since the following inequality
\[
\left|w+(x)-w^+(y)\right|\le |w(x)-w(y)|
\]
for all $(x,y)\in \R^N\times\R^N$ and function $w:\R^N\rightarrow\R$.
Therefore, $v_1\in X_0^\alpha(D_0)$. Similarly, we can show $v_2\in X_0^\alpha(D_0)$.
On the other hand, by Theorem 6 in \cite{FSV}, we know that $v_i$ can be approximate arbitrarily closely in the $X_0^\alpha(D_0)$ norm by $C_0^\infty(D_0)$ functions. Hence, we see that (\ref{2.2}) holds when $\phi_i$ is replaced by $v_i$ for  $i=1,2$.

Denote
\[
D(\ve)=\{x\in \Om: u_2(x)+\ve_2>u_1(x)+\ve_1\}.
\]
We notice that the integrands in the right hand side of (\ref{2.2}) (with $\phi_i=v_i$) vanishing outside $D(\ve)$.  Next, we prove the  left hand side of (\ref{2.2}) in nonpositive.  We first divide $\R^{2N}$ into four disjoint region as:
\begin{eqnarray*}
\R^{2N}&=&\left[\R^N\setminus D(\ve)\times\R^N\setminus D(\ve)\right]\cup\left[D(\ve)\times\R^N\setminus D(\ve)\right]\\
&\cup&\left[\R^N\setminus D(\ve)\times D(\ve)\right]\cup\left[D(\ve)\times D(\ve)\right].
\end{eqnarray*}
For $(x,y)\in \R^N\setminus D(\ve)\times\R^N\setminus D(\ve)$, we know that $v_i(x)=v_i(y)=0$, $i=1,2$. Therefore,
\[
A_1:=\int_{\R^N\setminus D(\ve)}\int_{\R^N\setminus D(\ve)}\frac{\mathcal{L}(u_1,u_2)}{|x-y|^{N+2\alpha}}dxdy=0,
\]
where 
\[
\mathcal{L}(u_1,u_2)=\left(u_1(x)-u_1(y)\right)\left(v_1(x)-v_1(y)\right)-\left(u_2(x)-u_2(y)\right)\left(v_2(x)-v_2(y)\right).
\]
For  $(x,y)\in D(\ve)\times\R^N\setminus D(\ve)$, we notice that $v_1(y)=v_2(y)=0$ and, by the definition of $D(\ve)$,
\beq\label{2.3}
u_2(x)+\ve_2> u_1(x)+\ve_1\quad{\rm and}\quad u_2(y)+\ve_2\le u_1(y)+\ve_1.
\eqq
It follows that
\begin{flushleft}
$\mathcal{L}(u_1,u_2)=\left[u_1(x)-u_1(y)\right]v_1(x)-\left[u_2(x)-u_2(y)\right]v_2(x)$\\
                                $=\left[(u_1(x)+\ve_1)-(u_1(y)+\ve_1)\right]v_1(x)-\left[(u_2(x)+\ve_2)-(u_2(y)+\ve_2)\right]v_2(x)$\\
                                $=\frac{[(u_2(x)+\ve_2)^2-(u_1(x)+\ve_1)^2]}{(u_1(x)+\ve_1)(u_2(x)+\ve_2)}\cdot[(u_1(x)+\ve_1)(u_2(y)+\ve_2)-(u_1(y)+\ve_1)(u_2(x)+\ve_2)]$\\
                                $\le0.$
 \end{flushleft}
Hence,
\[
A_2=\int_{D(\ve)}\int_{\R^N\setminus D(\ve)}\frac{\mathcal{L}(u_1,u_2)}{|x-y|^{N+2\alpha}}dydx\le0.
\]
A similar argument implies that
\[
A_3=\int_{R^N\setminus D(\ve)}\int_{ D(\ve)}\frac{\mathcal{L}(u_1,u_2)}{|x-y|^{N+2\alpha}}dydx\le0.
\]
Finally, if  $(x,y)\in D(\ve)\times D(\ve)$, it is easy to check that
\begin{eqnarray*}
\mathcal{L}(u_1,u_2)&=&\left(u_1(x)-u_1(y)\right)\left(v_1(x)-v_1(y)\right)-\left(u_2(x)-u_2(y)\right)\left(v_2(x)-v_2(y)\right)\\
&=&-(u_1(x)-u_1(y))^2+(u_1(y)-u_1(x))\left(\frac{(u_2(y)+\ve_2)^2}{u_1(y)+\ve_1}-\frac{(u_2(x)+\ve_2)^2}{u_1(x)+\ve_1}\right)\\
                                 &-&(u_2(x)-u_2(y))^2+(u_2(y)-u_2(x))\left(\frac{(u_1(y)+\ve_1)^2}{u_2(y)+\ve_2}-\frac{(u_1(x)+\ve_1)^2}{u_2(x)+\ve_2}\right).
\end{eqnarray*}
By Proposition \ref{p1}, we know that $\mathcal{L}(u_1,u_2)(x,y)\le0$ in $D(\ve)\times D(\ve)$. Therefore,
\[
A_4=\int_{D(\ve)}\int_{ D(\ve)}\frac{\mathcal{L}(u_1,u_2)}{|x-y|^{N+2\alpha}}dxdy\le0.
\]
Summing up these estimates from $A_1$ to $A_4$, we know that the left hand side of (\ref{2.2}) is nonpositive.

On the  other hand, as $\ve\rightarrow0$, the first term in the right hand side of (\ref{2.2}) converges to 
\[
\int_{D(0)} b(x)\left(u_2^{p-1}-u_1^{p-1}\right)(u_2^2-u_1^2)dx,
\]
while the last term in the right side of (\ref{2.2}) converges to 0. 

Next, we show that $D(0)=\emptyset$. Suppose  to the contrary that $D(0)\not=\emptyset$. Since the left side of (\ref{2.2})  is nonpositive by the estimates from $A_1$ to $A_4$ and right hand side of (\ref{2.2}) tends to 0 as $\ve\rightarrow0$, we easy deduce 
\[
\int_{\R^{2N}}\frac{\mathbb{L}(u_1,u_2)}{|x-y|^{N+2\alpha}}dxdy=0,
\]
where $\mathbb{L}(u_1,u_2)=\lim_{\ve\rightarrow 0}\mathcal{L}(u_1,u_2)$ and 
\[
\int_{D(0)} b(x)\left(u_2^{p-1}-u_1^{p-1}\right)(u_2^2-u_1^2)dx=0.
\]
This  imply that
\[b\equiv0\quad {\rm in }\,\, D(0)\]
and
\[
\mathbb{L}(u_1,u_2)\equiv 0\quad {\rm in }\,\, \R^N\times\R^N.
\]
Hence, by Proposition \ref{p1}, we know $u_1=ku_2$ in $D(0)$ for some constant $k$.
Since $b\not\equiv0$ in $\Om$, it follows from the above that $D(0)\not=\Om$. Thus, $D(0)\subset \Om$, $\partial D(0)\cap\Om\not=\emptyset$.
It follows that the open set $D(0)$ has connected component $\mathcal{G}$  such that $\partial \mathcal{G}\cap\Om\not=\emptyset$. Now on $\mathcal{G}$,
$u_1=ku_2$. On the other hand, we have $u_1|_{\partial \mathcal{G}\cap\Om}=u_2|_{\partial \mathcal{G}\cap\Om}>0$. Thus, $k=1$. So we have $u_1=u_2$ in 
$\mathcal{G}$, which contradicts $\mathcal{G}\subset D(0)$. Therefore, we must have $D(0)=\emptyset$
 and thus $u_1\ge u_2$ in $\Om$. We complete the proof of Lemma \ref{l2.1}.  $\Box$\\

By applying this comparison principle together with the Perron's method for the nonlocal equation, we can obtain the following two lemmas.

\begin{lemma}\label{l2.2}
Let $\Om$ be a bounded domain in $\R^N$ with smooth boundary and $p>1$. Suppose $a$ and $b$ are smooth positive functions in $\bar{\Om}$, and let $\mu_1$ denote the
first eigenvalue of $(-\Delta)^\alpha u=\mu a(x)u$ in $\Om$ with $u=0$ in $\R^N\setminus \Om$. Then equation
 \begin{eqnarray*}
\left\{\begin{array}{l@{\quad }l}
(-\Delta)^\alpha u=\mu u[a(x)-b(x)u^{p-1}]&{\rm in}\,\,\Om,\\
u=0&{\rm in}\,\,\R^N\setminus \Om
 \end{array}
 \right.
 \end{eqnarray*}
has a unique positive solution for every $\mu>\mu_1$. Furthermore,  the unique solution $u_\mu$ satisfies $u_\mu\rightarrow [a(x)/b(x)]^{1/{p-1}} $ 
uniformly in amy compact subset of $\Omega$ as $\mu\rightarrow +\infty$. 
\end{lemma}

{\bf Proof.}  (Existence) The existence follows from a simple sub- and super-solution argument. 
In fact, any constant great than  or equal to $M=\max_{\bar{\Om}}[a(x)/b(x)]^{1/(p-1)}$ is a super-solution.
Let $\phi$ be a positive eigenfunction corresponding to $\mu_1$ (for the existence of the first eigenvalue and corresponding eigenfunction has been obtained in  \cite{PQ} and \cite{SV}), then for each fixed $\mu>\mu_1$ and small positive $\ve$, $\ve\phi<M$ and is a
sub-solution. Therefore, by the  sub- and super-solution method (see \cite{S}), there exist at least one positive solution.

(Uniqueness) If $u_1$ and $u_2$ are two positive solutions, by Lemma \ref{l2.1}, we have $u_1\le u_2$  and $u_2\le u_1$ both hold in $\Om$.
Hence, $u_1=u_2$. This proves the uniqueness.

(Asymptotic behaviour) Given any compact subset $K$ of $\Om$ and any small $\ve>0$ such that $\ve< v_0(x)=[a(x)/b(x)]^{1/(p-1)}$ in $\Om$. Let 
\[
v_\ve(x)=\left\{
  \begin{array}{ll}
v_0(x)+\ve   &\rm{in}\,\, K,    \\
l(x)        &\rm{in}\,\, \Om\setminus K,    \\
0, &\rm{in}\,\,\R^N\setminus \Om,    
\end{array}
\right.
\]
where $l(x)$ is nonnegative function such that $v_\ve$ is smooth in $\Om$ and satisfying $D_0:=\rm{supp }(v_\ve)\subset\subset \Om$. Thus, for any $x\in \Om$,
\begin{eqnarray*}
|(-\Delta)^\alpha v_\ve(x)|&\le&\int_{\R^N}\frac{|v_\ve(x)-v_\ve(y)|}{|x-y|^{N+2\alpha}}dy\\
&=&\int_{\Om}\frac{|v_\ve(x)-v_\ve(y)|}{|x-y|^{N+2\alpha}}dy+\int_{\R^N\setminus\Om}\frac{|v_\ve(x)|}{|x-y|^{N+2\alpha}}dy\\
&\le& \int_{\Om}\frac{|v_\ve(x)-v_\ve(y)|}{|x-y|^{N+2\alpha}}dy+\left(\int_{\R^N\setminus\Om}\frac{1}{dist(y,\partial D_0)^{N+2\alpha}}dy\right)\|v_\ve\|_{L^\infty(\R^N)}\\
&\le& C
\end{eqnarray*}
 for some positive constant $C=C(\ve)$ since $v_\ve$ is smooth and $dist(\partial \Om, \partial D_0)\ge \gamma>0$.
 On the other hand,  we notice that $v_\ve(a(x)-b(x)v_\ve^{p-1})\le -\delta$ in $\Om$ for some positive constant $\delta=\delta(\ve)$.
  Hence, for all large $\mu$, $v_\ve$ is a super-solution
of our problem. 

On the other hand, let $\phi$ be a positive eigenfunction corresponding to $\mu_1$. Then we can find a small neighborhood
of $\partial \Om$ in $\Om$, say $U$, such that $\phi$ is very small in $U$. Therefore, for all $\mu>\mu_1+1$, we have 
\beq\label{u}
(-\Delta )^\alpha \phi=\mu_1 a(x)\phi\le \mu \phi[a(x)-b(x)\phi^{p-1}]\quad {\rm in}\,\, U.
\eqq
 By shrinking $U$ further if necessary, we can assume that $\bar{U}\cap K=\emptyset$ and $\phi<v_0-\ve$ in $U$. 
Next, we choose smooth function $w_\ve$ as 
\[
w_\ve(x)=\left\{
  \begin{array}{ll}
v_0(x)-\ve   &\rm{in}\,\, K,    \\
\phi(x)         &\rm{in}\,\, U,\\
l(x)        &\rm{in\,\, the \,\,rest\,\, of}\,\, \Om,    \\
0, &\rm{in}\,\,\R^N\setminus \Om,    
\end{array}
\right.
\]
where $l$ is a positive function such that $w_\ve$  is smooth in $\Om$ and satisfying $l\le v_0-\ve/2$. Moreover, we let 
\beq\label{phi}
\phi(x)\le w_\ve(x)\quad{\rm in}\,\,\Om
\eqq
otherwise we choose $\tilde{\phi}=\phi/C$ for some constant $C>0$ large replace $\phi$.
Then we can see that, for $x\in\Om\setminus U$, 
\begin{eqnarray*}
|(-\Delta)^\alpha w_\ve(x)|&\le&\int_{\R^N}\frac{|w_\ve(x)-w_\ve(y)|}{|x-y|^{N+2\alpha}}dy\\
&=&\int_{\Om}\frac{|w_\ve(x)-w_\ve(y)|}{|x-y|^{N+2\alpha}}dy+\int_{\R^N\setminus\Om}\frac{|w_\ve(x)-w_\ve(y)|}{|x-y|^{N+2\alpha}}dy\\
&=&\int_{\Om}\frac{|w_\ve(x)-w_\ve(y)|}{|x-y|^{N+2\alpha}}dy+\int_{\R^N\setminus\Om}\frac{|w_\ve(x)|}{|x-y|^{N+2\alpha}}dy\\
&=&\int_{\Om}\frac{|w_\ve(x)-w_\ve(y)|}{|x-y|^{N+2\alpha}}dy+\left(\int_{\R^N\setminus\Om}\frac{1}{dist(y,\partial U\cap\Om)^{N+2\alpha}}dy\right)\|w_\ve\|_{L^\infty(\R^N)}\\
&\le&C,
\end{eqnarray*}
 for some positive constant $C=C(\ve)$ since  $dist(\partial \Om, \partial U\cap \Om)\ge \gamma>0$. Moreover, we know $w_\ve(a(x)-b(x)w_\ve^{p-1})\ge \delta$ in $\Om\setminus U$ for some positive constant $\delta=\delta(\ve)$.
 Therefore,
 \beq\label{mu1}
 (-\Delta)^\alpha w_\ve \le \mu w_\ve(a(x)-b(x)w_\ve^{p-1})\quad{\rm in}\,\, \Om\setminus U
 \eqq
 for all large $\mu$. For $x\in U$, by (\ref{u}) and (\ref{phi}), we have
 \begin{eqnarray}\nonumber
 (-\Delta)^\alpha w_\ve (x)&=&\int_{\Om}\frac{w_\ve(x)-w_\ve(y)}{|x-y|^{N+2\alpha}}dy+\int_{\R^N\setminus\Om}\frac{w_\ve(x)}{|x-y|^{N+2\alpha}}dy\\\nonumber
 &=&\int_{\Om}\frac{\phi(x)-w_\ve(y)}{|x-y|^{N+2\alpha}}dy+\int_{\R^N\setminus\Om}\frac{\phi(x)}{|x-y|^{N+2\alpha}}dy\\\nonumber
  &\le&\int_{\Om}\frac{\phi(x)-\phi(y)}{|x-y|^{N+2\alpha}}dy+\int_{\R^N\setminus\Om}\frac{\phi(x)}{|x-y|^{N+2\alpha}}dy\\\nonumber
  &=& (-\Delta)^\alpha \phi(x)\\\nonumber
   &\le&\mu \phi[a(x)-b(x)\phi^{p-1}]\\\label{mu2}
   &=& \mu  w_\ve[a(x)-b(x)w_\ve^{p-1}],
\end{eqnarray}
for $\mu>\mu_1+1$. Finally, combining (\ref{mu1}) and (\ref{mu2}), we know $w_\ve$ is a sub-solution of our problem for all large $\mu$.

Since $w_\ve<v_\ve$, we deduce that $w_\ve\le u_\mu<v_\ve$ in $\Om$. In particular,
\[
[a(x)/b(x)]^{1/{(p-1)}} -\ve\le u_\mu\le [a(x)/b(x)]^{1/{(p-1)}} +\ve
\]  
in $K$ for all large $\mu$. Hence, $u_\mu\rightarrow [a(x)/b(x)]^{1/{p-1}} $ as $\mu\rightarrow +\infty$ in $K$, as required. $\Box$

\begin{lemma}\label{l2.3}
Let $\Om$, $a$ and $b$ be as in Lemma \ref{l2.2}.  Suppose $p$ verifies (\ref{p}), then equation
 \begin{eqnarray}\label{2.4}
\left\{\begin{array}{l@{\quad }l}
(-\Delta)^\alpha u=\mu u[a(x)-b(x)u^{p-1}]&{\rm in}\,\,\Om,\\
\lim_{x\in \Omega,x\rightarrow \partial\Om}u=+\infty,\\
u=g_\mu&{\rm in}\,\,\R^N\setminus \Om
 \end{array}
 \right.
 \end{eqnarray}
has at least one positive solution for each $\mu>0$ if the measurable function $g_\mu$ satisfying
\beq\label{g}
\int_{\R^N\setminus\Om}\frac{g_\mu(y)}{1+|y|^{N+2\alpha}}dy\le C,
\eqq
where positive constant $C$ is independent of $\mu$.
 Furthermore,  suppose $u_\mu$ is a positive solution of (\ref{2.4}), then $u_\mu$ satisfies $u_\mu\rightarrow [a(x)/b(x)]^{1/{(p-1)}} $ 
uniformly in amny compact subset of $\Omega$ as $\mu\rightarrow +\infty$. 
\end{lemma}

We first recall the following result in \cite{CFQ}.   Assume that $\delta>0$ such that the distance function $d(x)=dist(x,\partial \Om)$ is
of $C^2$ in $A_\delta=\{x\in \Om: d(x)<\delta\}$ and define
\beq\label{2.5}
V_\tau(x)=\left\{
  \begin{array}{ll}
l(x),   &x\in\Om\setminus A_\delta,    \\
d(x)^\tau, &x\in A_\delta,    \\
0,  &  x\in \R^N\setminus \Om, \\
\end{array}
\right.
\eqq
where $\tau$ is a parameter in $(-1,0)$ and the function $l$ is positive such that $V_\tau$ is $C^2$ in $\Om$.

\begin{proposition}\label{p2.2} (\cite{CFQ}, Proposition 3.2) 
Assume that $\Om$ is a bounded, open subset of $\R^N$ with a $C^2$ boundary. Then there exists $\delta_1\in (0,\delta)$ and s constant $C>1$ shch that
 if $\tau\in (-1,\tau_0(\alpha))$ where $\tau_0(\alpha)$ is the unique solution of
\[
C(\tau)=\int_0^{+\infty}\frac{\chi_{(0,1)}|1-t|^\tau+(1+t)^\tau-2}{t^{1+2\alpha}}dt
\]
for $\tau\in (-1,0)$  and $\chi_{(0,1)}$ is the characteristic function of the interval $(0,1)$,
 then 
\[
\frac{1}{C}d(x)^{\tau-2\alpha}\le -(-\Delta)^\alpha V_\tau(x)\le Cd(x)^{\tau-2\alpha},\quad for\,\, all \,\, x\in A_{\delta_1}.
\]
\end{proposition}

Next,  we will the existence result in Lemma \ref{l2.3} by applying Perrod's method and thus we need to find ordered sub and super-solution of (\ref{2.4}). As in \cite{CFQ}, we
begin with a simple lemma that reduce the problem to find them only in $A_\delta$. 

\begin{lemma}\label{l2.4} 
Let $\Om$, $a$ and $b$ be as in Lemma \ref{l2.2}.
Suppose $U$ and $W$ are order super and sub-solution of (\ref{2.4}) in the  sub-domain $A_\delta$. Then 
there exists $\lambda$ large such that $U_\lambda=U+\lambda \eta$ and $W_\lambda=W-\lambda \eta$ are ordered
super and sub-solution of (\ref{2.4}), where $\eta\in C_0^\infty(\R^N)$ satisfying  $0\le\eta\le1$ and $supp(\eta)\subset\Om\setminus A_\delta$.
\end{lemma}

{\bf Proof.} The proof is similar as Lemma 4.1in \cite{CFQ} and we just need replace $\bar{V}$ in  Lemma 4.1in \cite{CFQ} to $\eta$ for our lemma. So we omit the proof here. $\Box$\\

Now we in position to prove Lemma \ref{l2.3}. 

{\bf Proof of  Lemma \ref{l2.3}.} (Existence)  We define
\[
G_\mu(x)=\frac{1}{2}\int_{\R^N}\frac{\tilde{g}_\mu(x+y)}{|y|^{N+2\alpha}}dy \quad {\rm for}\,\, x\in\Om,
\]
where
\[
\tilde{g}_\mu(x)=\left\{
  \begin{array}{ll}
0,   &x\in\Om,    \\
g_\mu(x), & x\in \R^N\setminus \Om.\\
\end{array}
\right.
\]
We observe that 
\[
G(x)=-(-\Delta)^\alpha \tilde{g}_\mu(x)\quad {\rm for}\,\, x\in\Om.
\]
Moreover, we know that $G_\mu$ is continuous (see Lemma 2.1 in \cite{CFQ})  and nonnegative in $\Om$. Therefore, if $u$ is a solution of (\ref{2.4}), then $u-\tilde{g}_\mu$ is the solution of 
 \begin{eqnarray}\label{g}
\left\{\begin{array}{l@{\quad }l}
(-\Delta)^\alpha u=\mu u[a(x)-b(x)u^{p-1}]+G_\mu(x)&{\rm in}\,\,\Om,\\
\lim_{x\in \Omega,x\rightarrow \partial\Om}u=+\infty,\\
u=0&{\rm in}\,\,\R^N\setminus \Om
 \end{array}
 \right.
 \end{eqnarray}
and vice versa, if $u$ is a solution of (\ref{g}), then $u+\tilde{g}_\mu$ is a solution of (\ref{2.4}). Next, we will look for solution of (\ref{g}) instead of (\ref{2.4}).

Define 
\[
U_\lambda(x)=\lambda V_\tau(x)\quad {\rm and}\quad W_\lambda(x)=\lambda W_\tau(x),
\]
where $\tau=-2\alpha/(p-1)$. Notice that $\tau=-2\alpha/(p-1)\in (-1,\alpha-1)$, $\tau p=\tau-2\alpha$ and $\tau p<\tau<0$.

By Proposition \ref{p2.2}, we find that for $x\in A_\delta$ and $\delta>0$ small
\begin{eqnarray*}
(-\Delta)^\alpha U_\lambda&+&\mu b(x)U_\lambda^p-\mu a(x)U_\lambda-G_\mu(x)\\
&\ge& -C\lambda d(x)^{\tau-2\alpha}+\mu b(x)\lambda^p d(x)^{\tau p}-\mu a(x)\lambda d(x)^\tau-G_\mu(x)\\
&\ge & -C\lambda d(x)^{\tau-2\alpha}+\mu b(x)\lambda^p d(x)^{\tau p}-\mu a(x)\lambda d(x)^{\tau p}-G_\mu(x),
\end{eqnarray*}
for some $C>0$. Then there exists a large $\lambda>0$ such that $U_\lambda$ is a super-solution of (\ref{g}) with the first equation in $A_\delta$ since $G_\mu$ is continuous in $\Om$. Similarly, by Proposition \ref{p2.2}, we 
have that for $x\in A_\delta$ and $\delta>0$ small
\begin{eqnarray*}
(-\Delta)^\alpha W_\lambda&+&\mu b(x)W_\lambda^p-\mu a(x)W_\lambda-G_\mu(x)\\
&\le& -\frac{\lambda}{C} d(x)^{\tau-2\alpha}+\mu b(x)\lambda^p d(x)^{\tau p}-\mu a(x)\lambda d(x)^\tau-G_\mu(x)\\
&\le & -\frac{\lambda}{C} d(x)^{\tau-2\alpha}+\mu b(x)\lambda^p d(x)^{\tau p}\\
&\le& 0,
\end{eqnarray*}
if $\lambda>0$ small. Here we have used the fact $G_\mu$ is nonnegative.

Finally, by using Lemma \ref{l2.4}, there exists a solution $\tilde{u}_\mu$ of problem (\ref{g}) and thus a solution $u_\mu=\tilde{u}_\mu+\tilde{g}_\mu$ is a solution of (\ref{2.4}).
Moreover, $u_\mu>0$ in $\Om$. Indeed, since $0$ is a sub-solution of (\ref{2.4}), by Lemma \ref{l2.1}, we have $u_\mu\ge 0$ in $\Om$. If  $u_\mu(x^\prime)=0$ for some points  $x^\prime\in \Om$ and $u_\mu\not\equiv0$ in $\R^N$, then by the definition of fractional Laplacian $(-\Delta )^\alpha u_\mu(x^\prime)<0$ which is a contradiction. Therefore, $u_\mu>0$ in $\Om$.

(Asymptotic behaviour) 
Let $K$ be an arbitrary compact subset of $\Om$, $v_0(x)=[a(x)/b(x)]^{1/(p-1)}$ in $\Om$ and $\ve>0$ any small positive number satisfies $v_0>\ve$ in $\bar{\Om}$.
Define
\[
\tilde{w}_\ve(x)=\left\{
  \begin{array}{ll}
v_0(x)-\ve+\lambda\eta(x)  &{\rm in}\,\,K,    \\
\mu^{-1}d(x)^{\tau}              &{\rm in}\,\,A_\delta,    \\
l(x)                              &{\rm in\,\,the\,\,rest\,\,of}\,\,\Om,    \\
0 & {\rm in}\,\, \R^N\setminus \Om,\\
\end{array}
\right.
\]
where $\tau$ is a parameter in $(-1,0)$, $\lambda$ and $\eta$ defined as in Lemma \ref{l2.4} and  the function $l$ is positive such that $w_\ve$ is $C^2$ in $\Om$.
By a similar argument as Proposition 3.2 in \cite{CFQ}, there exists $\delta_1\in(0,\delta)$ and  constants $c>0$ and $C>0$ such that
\[
c(1+\mu^{-1}d(x)^{\tau-2\alpha})\le -(-\Delta)^\alpha \tilde{w}_\ve(x)\le C(1+\mu^{-1}d(x)^{\tau-2\alpha})
\]
for all  $x\in A_{\delta_1}$ and $\tau\in(-1,\alpha-1)$. Hence, for $x\in A_\delta$ and $\delta>0$ small,
\begin{eqnarray*}
(-\Delta)^\alpha \tilde{w}_\ve&+&\mu b(x)\tilde{w}_\ve^p-\mu a(x)\tilde{w}_\ve-G_\mu(x)\\
&\le& -c \mu^{-1}d(x)^{\tau-2\alpha}+\mu^{1-p} b(x)d(x)^{\tau p}- a(x)\lambda d(x)^\tau-G_\mu(x)-c\\
&\le & -c\mu^{-1}d(x)^{\tau-2\alpha}+\mu^{1-p} b(x)d(x)^{\tau p}\\
&\le& 0,
\end{eqnarray*}
if $\mu$ is large enough. Hence, $\tilde{w}_\ve$ is sub-solution in $A_\delta$. By applying Lemma \ref{l2.4}, we know that $w_\ve=\tilde{w}_\ve-\lambda\eta+\tilde{g}_\mu$ is a sub-solution of problem (\ref{2.4}) for all large $\mu>0$.

On the other hand, we define choose a  function 
\[
\tilde{v}_\ve(x)=\left\{
  \begin{array}{ll}
v_0(x)+\ve-\lambda\eta   &\rm{in}\,\, K,    \\
\mu d(x)^\tau      &\rm{in}\,\, A_\delta,\\
l(x)        &\rm{in\,\, the \,\,rest\,\, of}\,\, \Om,    \\
0, &\rm{in}\,\,\R^N\setminus \Om,    
\end{array}
\right.
\]
where $\tau$ is a parameter in $(-1,0)$, $\lambda$ and $\eta$ defined as in Lemma \ref{l2.4} and  the function $l$ is positive such that $v_\ve$ is $C^2$ in $\Om$.
By a similar argument as Proposition 3.2 in \cite{CFQ}, there exists $\delta_1\in(0,\delta)$ and  constants $c>0$ and $C>0$ such that
\[
c(1+\mu d(x)^{\tau-2\alpha})\le -(-\Delta)^\alpha \tilde{v}_\ve(x)\le C(1+\mu d(x)^{\tau-2\alpha})
\]
for all  $x\in A_{\delta_1}$ and $\tau\in(-1,\alpha-1)$. Hence, for $x\in A_\delta$ and $\delta>0$ small,
\begin{eqnarray*}
(-\Delta)^\alpha \tilde{v}_\ve&+&\mu b(x)\tilde{v}_\ve^p-\mu a(x)\tilde{v}_\ve-G_\mu(x)\\
&\ge& -C \mu d(x)^{\tau-2\alpha}+\mu^{p+1} b(x)d(x)^{\tau p}- \mu^2 a(x)\lambda d(x)^\tau-G_\mu(x)-C\\
&\ge & -C\mu d(x)^{\tau-2\alpha}+\mu^{p+1} b(x)d(x)^{\tau p}- \mu^2 a(x)\lambda d(x)^{\tau p}-G_\mu(x)-C\\
&\ge& 0,
\end{eqnarray*}
if $\mu$ is large enough since $\|G_\mu\|_{L^\infty(\bar{\Om})}\le C$ for some constant $C>0$ independent of $\mu$ by (\ref{g}). Hence, $\tilde{v}_\ve$ is sub-solution in $A_\delta$. By applying Lemma \ref{l2.4}, we know that $v_\ve=\tilde{v}_\ve+\lambda\eta+\tilde{g}_\mu$ is a super-solution of problem (\ref{2.4}) for all large $\mu>0$.

As $w_\ve<v_\ve$ in $\Om$, we must have $w_\ve\le u_\mu\le v_\ve$ in $\Om$.
This implies that $u_\mu\rightarrow v_0$ in $K$ as $\mu\rightarrow \infty$, as required. We complete the proof of this Lemma. $\Box$\\

 \setcounter{equation}{0}
\section{ Proofs}
 
 The main purpose of this section is to prove Theorems \ref{t1} and \ref{t2}. 
 
 \bigskip
 
{\bf Proof of Theorem \ref{t1}.}
Let us first observe that a nonnegative entire solution of (\ref{1.1}) is either identically zero or positive everywhere. Indeed, if  
$u(x^\prime)=0$ for some points  $x^\prime\in \R^N$ and $u\not\equiv0$ in $\R^N$, then by the definition of fractional Laplacian $(-\Delta )^\alpha u(x^\prime)<0$ which is a contradiction. 
Therefore, we only consider positive solution.

Suppose $\lambda>0$ and let $u$ be an arbitrary positive entire solution of (\ref{1.1}). We will show that $u(x_0)=\lambda^{1/(p-1)}$ for any point $x_0\in \R^N$ by using pointwise convergence 
of Lemmas \ref{l2.2} and \ref{l2.3}.

For any $t>0$, define
\[
u_t(x)=u[x_0+t(x-x_0)].
\]
Then $u_t$ satisfies
\[
(-\Delta)^\alpha u=t^{2\alpha}(\lambda u-u^p)\quad {\rm in}\,\, \R^N.
\]

Let $B$ denote the unit ball with center $x_0$. By Lemma \ref{l2.2}, for all large $t$, the problem
 \begin{eqnarray*}
\left\{\begin{array}{l@{\quad }l}
(-\Delta)^\alpha v=t^{2\alpha}u(\lambda-u^{p-1})&{\rm in}\,\,B,\\
u=0&{\rm in}\,\,\R^N\setminus B
 \end{array}
 \right.
 \end{eqnarray*}
 has a unique positive solution $v_t$ and $v_t\rightarrow \lambda^{1/(p-1)}$ as $t\rightarrow \infty$ at $x=x_0\in B$.
 By Lemma \ref{l2.1}, we have that $u_t\ge v_t$ in $B$ and thus
 \[
 u(x_0)=u_t(x_0)\ge v_t(x_0).
 \]
 Letting $t\rightarrow \infty$ in the above inequality we conclude that $u(x_0)\ge \lambda^{1/(p-1)}$.
 
 Let $w_t$ be a positive solution of 
  \begin{eqnarray*}
\left\{\begin{array}{l@{\quad }l}
(-\Delta)^\alpha w=t^{2\alpha} w(\lambda-w^{p-1})&{\rm in}\,\,B,\\
\lim_{x\in B,x\rightarrow \partial B}w=+\infty,\\
w=u_t&{\rm in}\,\,\R^N\setminus \bar{B}.
 \end{array}
 \right.
 \end{eqnarray*} 
By our assumption, we know that
 \beq\label{t}
 \int_{\R^N}\frac{u_t(x)}{1+|x|^{N+2\alpha}}dx\le C, 
 \eqq
 where constant $C>0$ independent of $t$ for $t$ large enough. In fact
 \begin{eqnarray}\nonumber
 \int_{\R^N}\frac{u_t(x)}{1+|x|^{N+2\alpha}}dx&=& \int_{\R^N}\frac{u(x_0+t(x-x_0))}{1+|x|^{N+2\alpha}}dx\\\label{tt}
 &=& \int_{\R^N}\frac{u(x)}{t^N\left(1+\left|\frac{x+(t-1)x_0}{t}\right|^{N+2\alpha}\right)}dx.
\end{eqnarray}
Define function
\[
f(t)=t^N\left(1+\left|\frac{x+(t-1)x_0}{t}\right|^{N+2\alpha}\right),
\]
we know $f(1)=1+|x|^{N+2\alpha}$ and $f(t)\rightarrow +\infty$ as $t\rightarrow+\infty$. Then, we can choose $t$ large enough such that $f(t)\ge f(1)$. 
So by (\ref{tt}), for $t$ large enough, we have
\[
\int_{\R^N}\frac{u_t(x)}{1+|x|^{N+2\alpha}}dx\le \int_{\R^N}\frac{u(x)}{1+|x|^{N+2\alpha}}dx\le C,
\]
since  $u\in L^1(\R^N,\omega)$.

 Then, applying Lemma \ref{l2.3}, we see that $w_t\rightarrow \lambda^{1/(p-1)}$ as $t\rightarrow \infty$ at $x=x_0\in B$.
 Applying Lemma \ref{l2.1}, we have that $u_t\le w_t$ in $B$ and thus
 \[
 u(x_0)=u_t(x_0)\le w_t(x_0).
 \]
 Letting $t\rightarrow \infty$ in the above inequality we conclude that $u(x_0)\le \lambda^{1/(p-1)}$. Therefore, $u(x_0)= \lambda^{1/(p-1)}$.
 Since $x_0$ is arbitrary, we conclude that $u\equiv \lambda^{1/(p-1)}$ in $\R^N$ for $\lambda>0$, the unique constant solution of (\ref{1.1}). $\Box$\\

 Next,  we will extend Theorem \ref{t1} to similar problem with variable coefficients, that is, Theorem \ref{t2}.
 We first consider the following equation which is more general than (\ref{1.2}): 
 \beq\label{3.3}
 (-\Delta)^\alpha u= a(x)u-b(x)u^p,\quad x\in \R^N, 
 \eqq
 where $a(x)$ and $b(x)$ are continuous functions in $\R^N$ and satisfying 
 \beq\label{3.4}
 \lim_{|x|\rightarrow \infty} a(x)=a^\infty>0, \quad  \lim_{|x|\rightarrow \infty} b(x)=b^\infty>0.
 \eqq
 Here we allow $a$ and $b$ can be change sign which is more general than (\ref{1.2}).
 
 \begin{teo}\label{t3.1}
 Under the above assumptions, if $u\in C_{loc}^{2\alpha+\beta}(\R^N)\cap L^1(\R^N,\omega)$ for some $\beta>0$ is a positive solution of (\ref{3.3})
 with $p$ verifies (\ref{p}), then 
 \[
 \lim_{|x|\rightarrow \infty}u(x)=\left(\frac{a^\infty}{b^\infty}\right)^{\frac{1}{p-1}}.
 \]
 \end{teo}
 
 We postpone the proof of Theorem \ref{t3.1}  and first we use it to prove the following result.
 
 \begin{corollary}\label{c3.1}
 Under the assumptions in Theorem \ref{t3.1}, if we further assume that $b$ is a nonnegative, then 
 problem (\ref{3.3}) has at most one positive solution.
 \end{corollary}
 
  {\bf Proof.} Suppose $u_1$ and $u_2$ are two positive solutions of (\ref{3.3}). By Theorem \ref{t3.1}, we have
  \[
 \lim_{|x|\rightarrow \infty}[(1+\ve)u_1-u_2]=\ve(a^\infty/b^\infty)^{1/(p-1)}>0
 \] 
 for any positive constant $\ve$. 
 
 Since $b$ is nonnegative, then $(1+\ve)u_1$ is a super solution of (\ref{3.3}). Therefore, applying Lemma \ref{l2.1} in a large ball to conclude
 that $(1+\ve)u_1\ge u_2$ in a large ball. It follows that this is true in all of $\R^N$. Hence, $u_1\ge u_2$ in $\R^N$ since $\ve$ is arbitrary. Similarly, we also can deduce
 $u_2\ge u_1$ in $\R^N$. Finally, we must have $u_1=u_2$  in $\R^N$, that is, (\ref{3.3}) has at most one positive solution. $\Box$\\
 
 Now we are in the position to prove Theorem \ref{t3.1}.
 
 {\bf Proof of Theorem \ref{t3.1}.}  We prove it by a contradiction argument. Assume that there exists a sequence points $x_n\in \R^N$ satisfying $|x_n|\rightarrow \infty$ such that
 $|u(x_n)-(a^\infty/b^\infty)^{1/(p-1)}|\ge \ve_0$ for some constant $\ve_0>0$.

 We define
 \[
 a_n(x)=a(x_n+x),\quad b_n(x)=b(x_n+x)\quad {\rm and}\,\, u_n(x)=u(x_n+x).
 \]
 Then $u_n$ satisfies
 \beq\label{3.5}
 (-\Delta)^\alpha u_n=a_n(x)u_n-b_n(x)u_n^p\quad{\rm in}\,\,\R^N.
 \eqq
 If we let 
 \[
 L_\alpha^nu=-(-\Delta)^\alpha u-(a^\infty-a_n)u,
 \]
 then we can rewrite (\ref{3.5}) as
 \[
 - L_\alpha^n u_n=a^\infty u_n-b_n(x)u_n^p\quad{\rm in}\,\,\R^N.
 \]

 Next, we fix a ball $B_r=\{x\in\R^N\,\,|\,\, |x|<r \}$ and consider the following problem
 \begin{eqnarray}\label{3.6}
\left\{\begin{array}{l@{\quad }l}
-L_\alpha^n w=a^\infty w-b_nw^p&{\rm in}\,\,B_r,\\
w=0&{\rm in}\,\,\R^N\setminus B_r.
 \end{array}
 \right.
 \end{eqnarray} 
 By using the variational characterization of the first eigenvalue and (\ref{3.4}), we see that $\lambda_1(-L^n_\alpha,B_r)\rightarrow \lambda_1((-\Delta)^\alpha,B_r)$ as $n\rightarrow \infty$, where $\lambda_1(-L^n_\alpha,B_r)$, $\lambda_1((-\Delta)^\alpha,B_r)$ denote the first eigenvalues of $-L^n_\alpha$ and $(-\Delta)^\alpha$ in $B_r$ with Dirichlet boundary conditions in $\R^N\setminus B_r$, respectively. Since we can choose $r$ large enough such that $ \lambda_1((-\Delta)^\alpha,B_r)<a^\infty$, then we may assume that $\lambda_1(-L^n_\alpha,B_r)<a^\infty$ for all $n$. On the other hand, we know that $b_n\rightarrow b^\infty$ uniformly in $B_r$ and thus we may also assume that $b_n\ge b^\infty/2$ in $B_r$ for all $n$.

 Let $\phi_n\in X_0^\alpha(B_r)$ be the first eigenfunction corresponding to $\lambda_1(-L_\alpha^n,B_r)$, that is,
  \begin{eqnarray}\label{phi}
\left\{\begin{array}{l@{\quad }l}
(-\Delta)^\alpha \phi_n+(a^\infty-a_n)\phi_n=\lambda_1(-L_\alpha^n,B_r)\phi_n&{\rm in}\,\,B_r,\\
\phi_n=0&{\rm in}\,\,\R^N\setminus B_r,
 \end{array}
 \right.
 \end{eqnarray} 
  with $\|\phi_n\|_{L^\infty(B_r)}=1$. By Theorems 1 and 2 in \cite{SV1} and using (\ref{3.4}),   we know $\phi_n$ is also a viscosity solution of (\ref{phi}).
 Then by Theorem 2.6 in \cite{CS1}, we have $\phi_n\in C^\beta_{loc}(B_r)$. Then, by Corollary 4.6 in \cite{CS}, $\phi_n$ converges uniformly to a $\phi_\infty$ and $\phi_\infty$ satisfies
    \begin{eqnarray*}
\left\{\begin{array}{l@{\quad }l}
(-\Delta)^\alpha \phi_\infty=\lambda_1((-\Delta)^\alpha,B_r)\phi_\infty&{\rm in}\,\,B_r,\\
\phi_\infty=0&{\rm in}\,\,\R^N\setminus B_r.
 \end{array}
 \right.
 \end{eqnarray*}  
  in viscosity sense. Next, by a similar argument as Theorem 2.1 in \cite{CFQ}, we know $\phi_\infty\in C^{2\alpha+\beta}_{loc}(B_r)$ and is a classical solution.
Then $\phi_\infty$ is the normalized positive eigenfunction corresponding to $\lambda_1((-\Delta)^\alpha,B_r)$.

 It is easily to check that $\ve\phi_n$ is a subsolution of (\ref{3.6}) for every $n$ if we choose $\ve$ small enough. Furthermore, $(2a^\infty/b^\infty)^{1/(p-1)}$ is a supersolution of (\ref{3.6})
 for all $n$. Then (\ref{3.6}) has a positive solution $w_n$ satisfies $\ve\phi_n\le w_n\le (2a^\infty/b^\infty)^{1/(p-1)}$. Then, using the regularity results again, we know $w_n$ converges in $C^{2\alpha+\beta}_{loc}(B_r)$  to some function $w$ satisfying $\ve\phi_\infty\le v\le (2a^\infty/b^\infty)^{1/(p-1)}$ and
    \begin{eqnarray*}
\left\{\begin{array}{l@{\quad }l}
(-\Delta)^\alpha w=a^\infty w-b^\infty w^p&{\rm in}\,\,B_r,\\
w=0&{\rm in}\,\,\R^N\setminus B_r.
 \end{array}
 \right.
 \end{eqnarray*}   
 Applying Lemma \ref{l2.2}, we know the above problem has a unique positive solution. Therefore, $w=w_r$ is uniquely determined and the whole sequence $w_n$ converges to $w_r$.

 By the comparison principle (see Lemma \ref{l2.1}), we know that
 \beq\label{3.7}
 u_n\ge w_n\rightarrow w_r\quad {\rm in}\,\, B_r .
 \eqq
 Next, we show $u_n$ has a uniformly bounded in $\R^N$  for all $n$ large enough, that is, there exists a positive constant $C$ independent of $n$ such that  $u_n(x_0)\le C$ for any $x_0\in \R^N$.
 We define, for any $t>0$, 
 \[
 u_{t,n}(x)=u_n[x_0+t(x-x_0)].
 \]
Then $u_{t,n}$ satisfying
\[
(-\Delta)^\alpha u=t^{2\alpha}(\tilde{a}_nu-\tilde{b}_nu^p)\quad{\rm in}\,\, \R^N,
\] 
 where $\tilde{a}_n(x)=a_n(x_0+t(x-x_0))$ and  $\tilde{b}_n(x)=b_n(x_0+t(x-x_0))$. On the other hand, since $\tilde{a}_n\rightarrow a^\infty$ and $\tilde{b}_n\rightarrow b^\infty$  uniformly 
 in $B$ where $B$ denote the unit ball with center $x_0$, we may assume $\tilde{a}_n\le 2a^\infty$ and $\tilde{b}_n\ge  b^\infty/2$ in $B$ for all $n$.

We consider the following problem
     \begin{eqnarray}\label{w1}
\left\{\begin{array}{l@{\quad }l}
(-\Delta)^\alpha v=t^{2\alpha}(2a^\infty v-(b^\infty/2)v^p)&{\rm in}\,\,B,\\
\lim_{x\in B, x\rightarrow\partial B}v=+\infty,\\
v=u_{t,n}&{\rm in}\,\,\R^N\setminus B.
 \end{array}
 \right.
 \end{eqnarray}   
 As a argument before, we know $u_{t,n}\in L^1(\R^N,\omega)$ for $t$ and $n$ large enough. Thus,
 by applying Lemma \ref{l2.3}, we know this problem has at least one positive solution. Let $v_t$ is a solution of (\ref{w1}), then $v_t\rightarrow (4a^\infty/b^\infty)^{1/(p-1)}$ as $t\rightarrow \infty$ at $x=x_0\in B$. Then the comparison principle  deduce that $u_{t,n}\le v_t$ in $B$ and thus
\[
u_n(x_0)=u_{t,n}(x_0)\le v_t(x_0).
\]
 Letting $t\rightarrow\infty$ in the above inequality we conclude that $u_n(x_0)\le (4a^\infty/b^\infty)^{1/(p-1)}$ as we required.

  Hence, $\|u_n\|_{L^\infty(\R^N)}\le C$ for all $n$ large enough, where constant $C>0$  independent of $n$.  On the other hand, $u_n\in C^{2\alpha+\beta}_{loc}(\R^N)$ implies that
  $u_n$ converges uniformly to some function $u_\infty$ and
  \[
(-\Delta)^\alpha u_n\rightarrow (-\Delta)^\alpha u_\infty\quad {\rm in}\,\,B_r
  \]
  is strongly as $n\rightarrow+\infty$.
  Hence, $u_\infty$ is nonnegative and  satisfies
 \[
 (-\Delta)^\alpha u=a^\infty u-b^\infty u^p\quad{\rm in }\,\, B_r.
  \]
 Furthermore, $u_\infty\ge w_r>0$. Thus $u_\infty$ is a positive solution and $|u_\infty(0)-(a^\infty/b^\infty)^{1/(p-1)}|\ge \ve_0$ due to the choice of $x_n$.

 Choose a sequence $r=r_1\le r_2\le \cdots \le r_m\rightarrow \infty$ as $m\rightarrow\infty$. We can apply the above argument to each $r_m$ and then use a diagonal
 process to obtain a positive solution $U$ of
 \beq\label{3.8}
 (-\Delta)^\alpha u=a^\infty u-b^\infty u^p\quad {\rm in}\,\, \R^N,
  \eqq
 which satisfies $U(0)\ge w_{r_m}(0)$ and $|U(0)-(a^\infty/b^\infty)^{1/(p-1)}|\ge \ve_0$. By changing of variables of the form $x=\theta y$, $\theta\in \R$, then (\ref{3.8}) can write
 as
  \[
 (-\Delta)^\alpha v=(a^\infty/b^\infty) v- v^p\quad {\rm in}\,\, \R^N,
  \]
 where $v(y)=u(x)=u(\theta y)$. In fact, we can choose $\theta=(b^\infty)^{-1/(2\alpha)}$. Then applying Theorem \ref{t1} to the above equation, we have $v\equiv (a^\infty/b^\infty)^{1/(p-1)}$. Hence, $ u\equiv (a^\infty/b^\infty)^{1/(p-1)}$.
 This a contradiction. We complete the proof. $\Box$\\

 {\bf Proof of Theorem \ref{t2}.} First, we let $\lambda>0$. We consider the following eigenvalue problem with weight function:
  \begin{eqnarray*}
\left\{\begin{array}{l@{\quad }l}
(-\Delta)^\alpha u=\lambda a(x)u&{\rm in}\,\,B_r,\\
u=0&{\rm in}\,\,\R^N\setminus B_r.
 \end{array}
 \right.
 \end{eqnarray*} 
 We denote $\mu_1$ be the first eigenvalue of this problem. Since $\mu_1\rightarrow 0$ as $r\rightarrow \infty$, we can choose $r_1>0$ large
 enough such that $\mu_1\le \lambda$ when $r\ge r_1$. So we can choose an increasing sequence $r_1<r_2<\cdots<r_n\rightarrow \infty$ and consider
 the following problem
  \begin{eqnarray}\label{3.1}
\left\{\begin{array}{l@{\quad }l}
(-\Delta)^\alpha u=\lambda a(x)u-b(x)u^p&{\rm in}\,\,B_n,\\
u=0&{\rm in}\,\,\R^N\setminus B_n,
 \end{array}
 \right.
 \end{eqnarray}  
 where $B_n=B_{r_n}$. By Lemma \ref{l2.2}, problem (\ref{3.1}) has a unique positive solution $u_n$ for each $n$.
 Furthermore, by the comparison principle (see Lemma \ref{l2.1}), we know $u_n\le u_{n+1}$. On the other hand, 
 any positive constant $M$ satisfying $M^{p-1}\ge M_0^{p-1}=\lambda\sup_{\R^N}a(x)/\inf_{\R^N}b(x)$ is a supersolution of (\ref{3.1}).
 It follows that $u_n\le M_0 $ for all $n$. Therefore, $u_n$ is increasing in $n$ and $u_\infty(x)=\lim_{n\rightarrow \infty}u_n(x)$ is well
 defined in $\R^N$. Then, $u_\infty$ satisfying (\ref{1.2}).
 Since $u_\infty\ge u_n>0$ in $B_n$ for each $n$, we know that $u_\infty$ is a positive solution of (\ref{1.2}). Moreover, by Corollary \ref{c3.1},
 $u_\infty$ is the unique solution of (\ref{1.2}). We complete the proof. $\Box$

 \setcounter{equation}{0}
\section{ Acknowledgements}
A. Quaas was partially supported by Fondecyt Grant No. 1151180 Programa Basal, CMM. U. de Chile and Millennium Nucleus Center for Analysis of PDE NC130017.

\end{document}